\documentclass[twoside,10pt]{amsart}
\usepackage{amsmath,amssymb}

\newcommand{\CP}{\mathbb{CP}}
\newcommand{\bCP}{\overline{\mathbb{CP}}}

\newcommand{\BZ}{\mathbb{Z}}

\newcommand{\TryPackage}[3]{\IfFileExists{#1.sty}{\usepackage{#1}#2}{#3}
}
\TryPackage{mathrsfs}{\renewcommand{\mathcal}{\mathscr}}{%
        \TryPackage{eucal}{}{}}

\newcommand{\ZZ}{{\mathbb Z}}
\newcommand{\RR}{{\mathbb R}}
\newcommand{\CC}{{\mathbb C}}

\usepackage{latexsym}
\usepackage{amsmath}
\usepackage{amsthm}
\usepackage{amscd}
\usepackage{xypic} 
\usepackage{amsfonts}
\usepackage{amssymb}
\usepackage{graphicx}
\usepackage{psfrag}

\newtheorem{df}{Definition}
\newtheorem{thm}[df]{Theorem}

\newtheorem{cor}[df]{Corollary}
\newtheorem{lem}[df]{Lemma}
\newtheorem{prop}[df]{Proposition}

\begin{document}

\title{Luttinger surgery and interesting  symplectic 4-manifolds  with small Euler characteristic}

\author{Scott Baldridge}
\author{Paul Kirk}
\date{January 23, 2007}

\thanks{The first   author  gratefully acknowledges support from the
NSF  grant DMS-0507857. The second  author  gratefully acknowledges
support from the NSF  grant DMS-0604310.}

\address{Department of Mathematics, Louisiana State University \newline
\hspace*{.375in} Baton Rouge, LA 70817}
\email{\rm{sbaldrid@math.lsu.edu}}

\address{Mathematics Department, Indiana University \newline
\hspace*{.375in} Bloomington, IN 47405}
\email{\rm{pkirk@indiana.edu}}


\maketitle

\begin{abstract}
In this article we construct a minimal symplectic 4-manifold $R$ that has small Euler characteristic ($e(R)=8$) and two essential Lagrangian tori with nice properties.  These properties make $R$ particularly suitable for constructing  interesting examples of symplectic manifolds with small Euler characteristic.  In particular, we construct an exotic symplectic $\CP^2\# 5\overline{\CP}^2$,  the smallest known minimal symplectic 4--manifold with $\pi_1=\BZ$,   the smallest known minimal symplectic 4-manifolds with $\pi_1=\BZ/a \oplus \BZ/b$ for all $a,b\in \BZ$, and the smallest known minimal symplectic 4-manifold with $\pi_1=\BZ^3$. We   use the  $\pi_1=\BZ$ example to derive a significantly better upper bound on the minimal Euler characteristic of all symplectic 4-manifolds with a prescribed fundamental group. 
\end{abstract}


\section{Introduction}

Combining the construction of taking  symplectic sums along a genus 2 surface in $(T^2\times S^2)\# 4\bCP^2$ (which we learned from  a recent  article  of   article of Akhmedov \cite{A} and also an article of Ozbaggi and Stipsitz \cite{OS}),   with the method of (symplectic) Luttinger surgery (\cite{Lut,ADK})  along Lagrangian tori,   in this article we prove the following theorem (see Theorem \ref{themanifold}).

\medskip
\noindent{\bf Theorem.} {\em There exits a minimal symplectic 4-manifold $R$ containing a pair of homologically essential Lagrangian tori  $T_1,T_2$. The manifold $R$ satisfies
\begin{itemize}
\item The Euler characteristic $e(R)=8$, the signature $\sigma(R)=-4$,
\item the meridians to $T_1$ and $T_2$ are nullhomotopic in $R-(T_1\cup T_2)$,
\item  $\pi_1(R)=\ZZ^2$ and the map $\pi_1(R-(T_1\cup T_2))\to \pi_1(R)$ is an isomorphism,
\item the inclusions induce the homomorphisms $\ZZ^2=\pi_1(T_1)\to \pi_1(R), (s,t)\mapsto (s,0)$ and
$\ZZ^2=\pi_1(T_2)\to \pi_1(R), (s,t)\mapsto (0,t)$.
\end{itemize}}

\medskip
  Combined with techniques such as Gompf's symplectic sum \cite{Gompf} and Luttinger surgery \cite{ADK} we use $R$ to construct interesting examples of minimal symplectic 4-manifolds with prescribed fundamental group and small Euler  characteristic.     In particular, we construct:
 \begin{enumerate}
\item A symplectic manifold $M$ with $\pi_1(M)=\ZZ$, $e(M)=8$, $\sigma(M)=-4$. This manifold contains a symplectic torus $T$ with trivial normal bundle which can be used as a smaller replacement to the elliptic surface $E(1)$ to kill   elements in fundamental groups of symplectic 4-manifolds  (Theorem \ref{coolthm}).

\item   For each group $G=\ZZ/a\oplus \ZZ/b$, with $a,b$ arbitrary integers, an infinite  family of minimal symplectic 4-manifolds with fundamental group $G$, Euler characteristic $8$ and signature $-4$ (Theorem \ref{cyclic}).  The case $a=b=1$ yields minimal symplectic 4-manifolds homeomorphic to 
$\CP^2\# 5\bCP^2$.

\item For each group $G$ presented with $g$ generators and $r$ relations, a minimal symplectic 4-manifold with fundamental group $G$, Euler characteristic $12+8(g+r)$ and signature $-8-4(g+r)$. This significantly improves the main result of \cite{BK}.
\end{enumerate}

The authors would like to thank A. Akhmedov, R. Fintushel, C. Judge,  C. Livingston, and J. Yazinsky for helpful discussions.

\section{ Some topology}

Consider the union of three symplectic surfaces $(T^2\times \{s_1\})\cup (\{r\} \times S^2)\cup (T^2\times \{s_2\})$  in the symplecctic 4-manifold $T^2\times S^2$ (with the product symplectic form). Resolving the three double points symplectically (\cite{Gompf}) yields a genus 2 symplectic surface $F$ with self-intersection 4. Blowing up $T^2\times S^2$ along four points on $F$  and taking the proper transform yields a genus 2 surface (which we continue to call $F$) in 
$$W=(T^2\times S^2)\# 4 \bCP^2.$$

The symplectic surface $F\subset W$ has the easily verified properties:
\begin{enumerate}
\item Based loops $a_1,b_1, a_2,b_2$  on $F$ representing the standard symplectic generators can be chosen so that $\pi_1(F)= \langle a_1, b_1, a_2, b_2\ | \ [a_1,b_1][a_2,b_2]=1 \rangle$ in such a way that $\pi_1(F)$  surjects to $\pi_1(W)=\ZZ a\oplus \ZZ b$    by the assignment
\begin{equation*}  a_1\mapsto a, b_1\mapsto b, a_2\mapsto a^{-1}, b_2\mapsto b^{-1}.
\end{equation*}
 \item $F$ intersects each of the four $-1$ exceptional spheres transversely once.
 
 \item $F$ has a trivial nomal bundle.
\end{enumerate}

Fix a trivialization of the normal bundle of $F$, and hence an identification of the boundary of a tubular neighborhood of $F$ with $F\times S^1$.   The meridian of $F$ is the curve $\{p\}\times S^1=\mu_F\subset W-F$ which is the boundary of a small normal disc  to $F$.  Up to a (free) homotopy, we may assume that $\mu_F$ lies on one of the exceptional spheres, and  hence $\mu_F$ is nullhomotopic in $W-F$. In particular, since $F$ is connected, every homotopy that  intersects $F$ can be replaced by a homotopy that misses $F$.  Since   every loop in $W$ can be pushed off $F$  it follows that the inclusion  $W-F\subset W$ induces an isomorphism  on fundamental groups.  Hence $\pi_1(W-F)=\ZZ a \oplus \ZZ b$, and the push off of $F$ into $W-F$ using the trivialization again induces the homomorphism 
$a_1\mapsto a, b_1\mapsto b, a_2\mapsto a^{-1}, b_2\mapsto b^{-1}.$

(Matsumoto \cite{M} describes a   Lefschetz fibration $f:W\to S^2$  with generic fiber $F$  and with eight singular fibers; the monodromy given by  the relation $(D_1D_2D_3D_4)^2=1$ in the mapping class group of  $F$, where $D_i$ is the positive Dehn twist about the curve $C_i$ of \cite[pg 325]{GS}.) 

\begin{lem}\label{lem1}
Let $R$ be any 4-manifold containing a genus 2 surface $G$ with trivialized normal bundle. Let $\phi:F\to G$
be a diffeomorphism, and set  $g_i=\phi_*(a_i),h_i=\phi_*(b_i)$ in $\pi_1(R)$.  Given a map  $\tau:F\to S^1$, let  $\tilde{\phi}:F\times S^1\to G\times S^1$ the diffeomorphism given by $\tilde{\phi}(a,s)=(\phi(a), \tau(a)\cdot s)$.
 Then the sum
$$S=R\#_{F,G}W= (R-nbd(G))\cup_{\tilde{\phi}}(W-nbd(F)))$$ has fundamental group
$$\pi_1(S)=\pi_1(R)/N(g_2g_1, h_2h_1, [g_1,h_1])$$
where $N(g_1g_2, h_1h_2,[g_1,h_1])$ denotes the normal subgroup generated by $g_1g_2$, $h_1h_2,$ and the commutator $[g_1,h_1]$.
\end{lem}

\begin{proof}

The homomorphism $\pi_1(R-G)\to \pi_1(R)$ is  a surjection since every loop can be pushed off $G$, and the kernel is normally generated by the meridian $\mu_G=\tilde{\phi}(\mu_F)$ because any homotopy can be made transverse to $G$ and $G$ is connected.

 The Seifert-Van Kampen theorem   implies that $\pi_1(S)$ is the quotient of the free product
 $\pi_1(R-G)*\pi_1(W-F)= \pi_1(R-G)*(\ZZ a\oplus \ZZ b)$ by the normal subgroup generated by $ \tilde{\phi}(\mu_F)\mu_G^{-1},  g_1 a^{-1},   h_1 b^{-1}, g_2a , h_2b $. Since $\mu_F=1$, this  can also be described as the quotient of $\pi_1(R )*(\ZZ a\oplus \ZZ b)$ by the normal subgroup generated by
$ g_1 a^{-1},   h_1 b^{-1}, g_2a , h_2b$. We can then eliminate the generators $a$ and $b$ and conclude that $\pi_1(S)$ is the quotient of $\pi_1(R)$ by the normal subgroup generated by $g_2g_1, h_2h_1,$ and the commutator $[g_1,h_1]$.

\end{proof}

We next construct a useful building block for our subsequent constructions. We will postpone the discussion of symplectic structures until the next section.  Our emphasis now is on a careful calculation of fundamental groups.

We begin with the statement of a well-known lemma which computes the fundamental groups of generalized mapping tori. We omit the standard proof.

\begin{lem}\label{HNN} Let $X$ be a finite connected CW complex with base point $x_0$, $i:Y\subset X$  the inclusion of
a connected subcomplex containing $x_0$, and $f:Y\to X$ a continuous base point preserving map. Suppose $\pi_1(X,x_0) $ is generated by loops $x_1,\cdots, x_g$, with a presentation $\langle x_1,\cdots, x_g\ | \ w_1,\cdots , w_r\rangle$, and $\pi_1(y,x_0)$ is generated by loops $\gamma_1,\cdots, \gamma_n$.

Let $Z=X\times [0,1]/\sim$ where $(y,0)\sim(f(y),1)$ for $y\in Y$. Give $Z$ the base point $z_0=(x_0,0)$ and denote by $x_i\in \pi_1(Z,z_0)$ again the images of the loops $x_i$ under the map
$X=X\times\{0\}\subset X\times [0,1]\to Z$. Let $t\in \pi_1(Z,z_0)$ denote the loop
$u\mapsto (x_0,u)\subset X\times [0,1]\to Z$.

Then $\pi_1(Z,z_0)$ has the presentation
$$\pi_1(Z,z_0)=\langle t,  x_1,\cdots,  x_g\ | \ w_1,\cdots, w_r, t i_*(\gamma_1)t^{-1}=f_*(\gamma_1),
\cdots,  ti_*(\gamma_n)t^{-1}=f_*(\gamma_n)\rangle.$$
\qed
\end{lem}

Let $H$ be a once punctured torus. Choose a base point $h\in H$ on the boundary of $H$  and denote by $x,y\in \pi_1(H,h)$ the two standard symplectic free generators of the fundamental group. Let $D:H\to H$ denote the Dehn twist along a curve parallel to $x$.  Let $K\subset H$ denote a loop parallel to $y$ and $L$ a loop parallel to $x$,  as illustrated in the following figure.

\begin{figure}[h]
\begin{center}
\psfrag{y}{$y$} \psfrag{x}{$x$}
\psfrag{D}{$D$} \psfrag{h}{$h$}\psfrag{K}{$K$}\psfrag{L}{$L$}
   \includegraphics{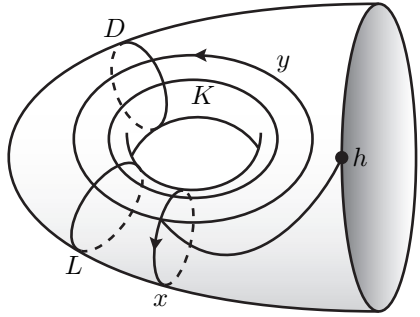}
   \caption{The surface $H$.}
   \label{fig2}
\end{center}
\end{figure}

Let $Z$ denote the mapping torus of $D:H\to H$, i.e.
$$Z=H\times[0,1]/\sim \text{ where } (x,0)\sim (D(x),1).$$
Thus $Z$ is a fiber bundle over $S^1=[0,1]/(0\sim 1)$ with fibers $H$.  Let $C=Z\times S^1$. For convenience we think of the second coordinate as $[0,1]/(0\sim 1)$. Hence we give $C$ coordinates $(a,u,v)\in H\times [0,1]\times [0,1]$.   We have loops $x=x\times(0,0), y=y\times(0,0), t=\{h\}\times I\times \{0\}$, and
$s=\{h\}\times \{0\}\times I$ all based at $h=(h,0,0)$ (note that $D$ fixes $h$, so that $t$ and $s$ are loops). Applying Lemma \ref{HNN} gives:

$$\pi_1(C)=\langle x,y,t\ | txt^{-1}=x, tyt^{-1}= yx\rangle\oplus \ZZ s=
\langle x,y,t\ | \ [t,x], [y^{-1},t]x^{-1}\rangle\oplus \ZZ s.$$

\medskip
There are two interesting tori in $C$. One is the torus $T_1=K\times\{0\}\times I/\sim,$ the other is the torus $T_2=L\times I\times \{0\}/\sim$.  Notice that the definition of $T_2$ makes sense since $L$ misses the support of the Dehn twist $D$.

We connect $T_1$ and $T_2$ to the base point $(h,0,0)$ as follows. For $T_1$, follow the path in $H\times\{(0,0)\}$ starting at $(h,0,0)$ illustrated in Figure 1 to the intersection of the loops labeled $x$ and $y$, then follow $x$ in the reverse direction a short way until you reach $K=K\times\{(0,0)\}\subset T_1$.   For $T_2$  follow the path to the intersection of $x$ and $y$, then follow $y$ backwards until you hit $L= L\times \{(0,0)\}\subset T_2$.
Based this way, the two generators of $\pi_1(T_1)$ map to the classes $y$ and $s$ in $\pi_1(C)$, and the two generators of $\pi_1(T_2)$ map to the classes $x$ and $t$ in $\pi_1(C)$.

We are interested in the fundamental group of their complement and corresponding meridian circles. We begin with the calculation of
$\pi_1(C-T_1)$. Notice that (identifying $K$ with $K\times \{0\}$)
 $C-T_1=(Z-K)\times S^1$, and hence $\pi_1(C-T_1)=\pi_1(Z-K)\oplus \ZZ s$. Moreover the meridian  $\mu_1$ of  $T_1$ (the homotopy class of the boundary of a small 2-disk transverse in $C$ to $T_1$) is represented by a loop in $Z-K$.

 Since $Z-K$ is constructed as the quotient of $H\times [0,1]$ where $(x,0)\sim (D(x),1)$  for $x\in H-K$,
 and $\pi_1(H-K)$ is generated by $y$ and $c=[x,y]=\partial H$, it follows  from Lemma \ref{HNN} that
 $$\pi_1(C-T_1)=  \langle x,y,t\ | \  [y^{-1},t ]x^{-1}, [t,[x,y]]\rangle \oplus \ZZ s$$
 and the map $\pi_1(C-T_1)\to \pi_1(C)$ takes   $x$ to $x$, $y$ to $y$, $t$ to $t$, and $s$ to $s$.
The meridian $\mu_1$ is represented  by $[t,x]$.

A similar argument computes $\pi_1(C-T_2)$: this time
$C-T_2= Z\times [0,1]/\sim$, where $(z,0)\sim(z,1)$ when $z\in Z- T_2$. Thus the generator $s$ only commutes with those elements of $\pi_1(Z)$ which are represented by loops that miss $T_2$.  Lemma \ref{HNN} implies
$$\pi_1(C-T_2)
=\langle x,y,t,s\ | \ [t,x], [y^{-1},t]x^{-1}, [s,x], [s,t]\rangle $$
(with the obvious morphism to $\pi_1(C)$) and the meridian $\mu_2$ is represented by $[s,y]$.

Finally, note that the boundary $\partial C$ is a 3-torus, with $\pi_1(\partial C,h)=\ZZ c \oplus \ZZ t\oplus \ZZ s$ with $c=[x,y]$.

\medskip

\noindent{\bf Remark.} The only subtle parts of these calculations are the computations of $\mu_1$ and $\mu_2$. However, an alternative and simpler calculation which is quite sufficient for our purposes is the following. Notice that the torus $T_1$ intersects $T_2$ transversally in one point, say $p$. Choose a pair of loops on $T_2$ which miss $p$  and are freely homotopic to $x$ and $t$. Connect these loops to the base point of $C$ in $C-T_1$, and call the result $x'$ and $t'$. Then obviously $[t',x']$ is a meridian for $T_1$.  Moreover, $t'$ is conjugate to $t$ and $x'$ is conjugate to $x$.  To anticipate what follows, notice that if $C-T_1$ is modified in some way so that $x$ becomes nullhomotopic, then $x'$ also becomes nullhomotopic and hence also $[t',x']$.  The same argument shows that a meridian for $T_2$ is of the form $[s',y']$ for some curves $s'$ conjugate to $s$ and $y'$ conjugate to $y$.

\section{The symplectic manifold $R$}\label{sec3}

Let $G$ denote another surface of genus $2$, with base point $g$.     Fix based loops $x_1,y_1,x_2,y_2$ representing a symplectic basis of $\pi_1(G)$. Let $D_1:G\to G$ denote the Dehn twist about a curve parallel to $x_1$, and let $D_2:G\to G$ denote the Dehn twist about a curve parallel to $y_2$.  Let $K$ be a curve parallel to $y_1$ and $L$ a curve parallel to $y_2$. These curves are illustrated in the following figure.

\begin{figure}[h]
\begin{center}
\psfrag{y1}{$y_1$} \psfrag{x1}{$x_1$}\psfrag{y2}{$y_2$}
\psfrag{x2}{$x_2$} \psfrag{D1}{$D_1$}\psfrag{D2}{$D_2$}\psfrag{g}{$g$}\psfrag{K}{$K$}\psfrag{L}{$L$}
   \includegraphics{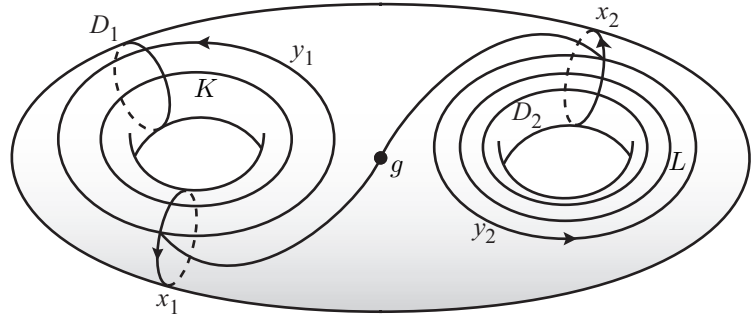}
   \caption{The surface $G$. }
   \label{fig:surface_genus_g}
\end{center}
\end{figure}

Denote by $\phi:G\to G$ the composite $\phi=D_2\circ D_1$.   The mapping torus of $\phi$ is a 3-manifold   $$Y=G\times [0,1]/\sim, \  \text{ where } (x,0)\sim (\phi(x),1) \text{ for } x\in G.$$ The projection $Y\times [0,1]\to [0,1]$ determines a fibration $p:Y\to S^1$ of $Y$ over the circle.    Let $B=Y\times S^1$.

Notice that $B$ is the union along their boundary  $T^3$ of two copies of the manifold $C$ constructed in the previous section, but with a reversal of the roles of $x$ and $y$ in the second copy (to the right in the figure). The Seifert-Van Kampen theorem  coupled with the calculation of the previous section give
$$\pi_1(B)=\langle x_1,y_1,x_2,y_2, s,t \ | \  [x_1,y_1][x_2,y_2], [t,x_1], [y_1^{-1},t]x_1^{-1},
[x_2^{-1},t]y_2^{-1}, [t,y_2]\rangle\oplus \ZZ s.$$

As before, we view $B$ as a quotient of $G\times [0,1]\times [0,1]$. We then have three disjoint surfaces in $B$: a pair of tori
$$T_1=K\times \{0\}\times I/\sim$$
$$T_2=L\times I\times \{0\}/\sim$$
and a genus two surface
$$G=G\times \{\tfrac12\}\times  \{\tfrac12\}/\sim.$$

The tori $T_1$ and $T_2$ can be connected to the base point $(g,0,0)$ by the chosen  path in $G\times\{(0,0)\} $, and  $G$ can be connected to the base point using the path $(g,u,u),u\in[0,\tfrac12]$.
The $T_i$  correspond to the tori identified in the previous section (again with an appropriate reversal of $x$ and $y$ in the second copy of $C$), and so the Seifert-Van Kampen theorem  can be used to compute $\pi_1(B-(T_1\cup T_2))$ and the two corresponding meridians $\mu_1$, $\mu_2$.

 We summarize the resulting fundamental group calculations in the next   proposition.

 \begin{prop}\label{prop1} For the manifold $B$ constructed above and the three surfaces $T_1,T_2$ and $G$ we have
 $$\pi_1(B)=\langle x_1,y_1,x_2,y_2, t \ | \  [x_1,y_1][x_2,y_2], [t,x_1], [y_1^{-1},t]x_1^{-1},
[x_2^{-1},t]y_2^{-1}, [t,y_2]\rangle\oplus \ZZ s,$$
$\pi_1(G)$ is generated by $x_1,y_1,x_2,$ and $y_2$, $\pi_1(T_1)$ is generated by $y_1$ and $s$, $\pi_1(T_2)$ is generated by $y_2$ and $t$.

Moreover, $\pi_1(B-(T_1\cup T_2))$ is generated by $  x_1,y_1,x_2,y_2, t$ and $s$ subject to the relations
$$ [x_1,y_1][x_2,y_2]=  [y_1^{-1},t]x_1^{-1}=
[t,[x_1,y_1]]= [x_2^{-1},t]y_2^{-1}= [t,y_2]=[s,x_1]=[s,y_1]= [s,t]=[s,y_2]=1.$$
In this group the meridian of $T_1$ is   (up to conjugation and change of orientation) $\mu_1=[x_1,t]$ and  the meridian of $T_2$ is $\mu_2=[x_2,s]$.
\qed
 \end{prop}

\medskip

There is a standard procedure, due to Thurston \cite{T}, for constructing a symplectic structure on $B$ for which   $G $   is a symplectic submanifold, and for which the tori $T_1$ and $T_2$ are Lagrangian.  Since it is critical for us that these tori be Lagrangian,   we give a proof of this fact for the convenience of the reader.

\begin{lem}\label{lem4} Let $A$ denote an annulus $[0,1]\times S^1$. Then there exists an area form $\alpha$ on $A$ and a diffeomorphism $D:A\to A$ isotopic to a Dehn twist supported away from a neighborhood of the boundary   so that $D^*(\alpha)=\alpha$.
\end{lem}
\begin{proof} Let $g:[0,1]\to [0,2\pi]$ be a smooth non-decreasing function so that $g(x)=0$  for $x<\tfrac{1}{3}$ and $g(x)=2\pi$ for $x>\tfrac23$. Define $\tilde{D}:[0,1]\times\RR\to [0,1]\times\RR$ by $\tilde{D}(x,y)=(x,y+g(x))$. Then $\tilde{D}^*(dx\wedge dy)=dx\wedge dy $. Since $\tilde{D}(x,y+2\pi)=\tilde{D}(x,y)$,
$\tilde{D}$ descends to a diffeomorphism of $[0,1]\times S^1$ which   represents a Dehn twist, preserves the standard area form, and is supported away from the boundary.
\end{proof}

 \begin{cor} There exists a symplectic structure on $B$ so that $G$ is symplectic and $T_1$ and $T_2$ are Lagrangian.  Moreover, $B$ is minimal.
 \end{cor}
 \begin{proof} By Lemma \ref{lem4} we can find a Riemannian metric $g$ on the surface $G$ so that  the Dehn twists $D_1$ and $D_2$ preserve the area form $\alpha(g)$ and are supported in small annular neighborhoods that miss the curve   $L$.   Let $q_i$ $i=1,2,3$ denote the projections of
 $G\times[0,1]\times[0,1]$.   Then $\omega= q_1^*(\alpha(g)) + q_2^*(dt)\wedge q_3^*(ds)$ is a symplectic form
 on $G\times[0,1]\times[0,1]$ for which the slices $G\times\{(t,s)\}$ are symplectic and $K\times\{0\}\times I$ and $L\times I\times\{0\}$ are Lagrangian. Since $(D_2D_1)^*(\alpha(g))=\alpha(g)$ and since $dt$ and $ds$ descend to $S^1$, $\omega$ descends to a well defined symplectic form on $B$ for which $G$ is symplectic and $T_1$ and $T_2$ are Lagrangian.

 To see that $B$ is minimal, it suffices to observe that the universal cover of $B$ is contractible, and so  $\pi_2(B)=0$, since this implies that there are no  spheres of self-intersection $-1$.
  \end{proof}

Rescale the symplectic form on Matsumoto's manifold $W$ if necessary so that the symplectic fiber $F=F_p  \subset W$ is symplectomorphic to the symplectic surface $G \subset B$. Fix a symplectomorphism which  takes the generators $a_1,b_1,a_2,b_2$ of $\pi_1(F)$ to $x_1,y_1,x_2,y_2\in \pi_1(G)$ respectively.

Then form the symplectic sum $R=W\#_{F,G}B$ of $W$ and $B$ along $F$ and $G$. As explained in Lemma \ref{lem1}, the fundamental group of $R$ is obtained from $\pi_1(B)$ by setting $x_1=x_2^{-1}, y_1=y_2^{-1}$, and requiring $x_1$ and $y_1$ to commute. Thus:
 \begin{eqnarray*}
 \pi_1(R)&=&\langle x_1,y_1,t\ | \ [x_1,y_1], [x_1,y_1][x_1^{-1},y_1^{-1}], [t,x_1], [y_1^{-1},t]x_1^{-1}, [x_1,t]y_1, [t,y_1^{-1}]\rangle\oplus \ZZ s\\
 &=&\ZZ t \oplus \ZZ s.
\end{eqnarray*}
In fact, the fourth and sixth relations force $x_1=1$, and then the fifth relation forces $y_1=1$.

The complement of $T_1\cup T_2$ in $R$ is the sum of $W$ with $B-(T_1\cup T_2)$ along $F$ and $G$.  Taking the presentation of $\pi_1(B-(T_1\cup T_2))$ of Proposition \ref{prop1},   setting $x_1=x_2^{-1}, y_1=y_2^{-1}$, and requiring $x_1$ and $y_1$ to commute, one sees  that $[t,y_1^{-1}]=1$ and $[y_1^{-1},t]x_1^{-1}=1$. This forces $x_1=1$. Then $[x_1,t]y_1=1$ forces $y_1=1$. The presentation thus reduces to
$$
 \pi_1(R-(T_1\cup T_2))=\ZZ t \oplus \ZZ s.
$$
 In particular $\mu_1$ and $\mu_2$ are trivial in $\pi_1(R-(T_1\cup T_2))$.

Gompf's symplectic sum theorem \cite{Gompf} guarantees that a symplectic structure can be found on $R$ for which the tori $T_1$  and $T_2$ remain Lagrangian.
Notice that each $T_i$ represents a non-zero class in $H_2(R;\RR)$. In fact, the union of a small $D^2$ transverse to  $T_1$ and a nullhomotopy of $\mu_1$ in $R-(T_1\cup T_2)$  gives a (singular)
2-sphere intersecting $T_1$ transversally in one point. This shows that $T_1$ represents a non-zero class in $H_2(R;\RR)$ and $H_2(R-T_2;\RR)$. Similarly $T_2$ represents a non-zero class in $H_2(R;\RR)$ and $H_2(R-T_1;\RR)$.  The classes $[T_1]$ and $[T_2]$ are not multiples of each other since the singular 2-sphere  dual to $T_1$ misses $T_2$.  It follows \cite[Lemma 1.6]{Gompf} that  the symplectic form on $R$ can be perturbed by an arbitrarily small 2-form so that $T_1$ and $T_2$ become symplectic.

\medskip

The Euler characteristic  of $R$ is computed using the formula $e(A\#_HB)=e(A)+e(B)-2e(H)$ valid for any symplectic  sum of 4-manifolds.  From this formula we first see that $e(R)=e(W)+e(B)+4=4+0+4=8$.
Novikov additivity computes the signature $\sigma(R)=\sigma(W)+\sigma(B)=-4$.
We will show in Theorem \ref{theyaremini}  below that $R$ is minimal.

\medskip

We summarize the properties of the manifold $R$ in the following theorem

\begin{thm}\label{themanifold}  There exists a closed minimal symplectic 4-manifold $R$  with Euler characteristic  $e(R)=8$ and signature $\sigma(R)=-4$ containing a disjoint  pair of Lagrangian, homologically essential  tori $T_1$, $T_2$ with  trivial normal bundles so that $\pi_1(R)=\ZZ t\oplus \ZZ s$ and such that the homomorphism induced by the inclusion $\pi_1(R-(T_1\cup T_2))\to \pi_1(R)$ is an isomorphism.

Moreover, the meridians $\mu_1$ and $\mu_2$ to the tori $T_1$ and $T_2$ are nullhomotopic in $R-(T_1\cup T_2)$. The homomorphism  induced by inclusion  $\pi_1(T_1)\to \pi_1(R)$ takes one generator to s and the other to the identity and the homomorphism $\pi_1(T_2)\to \pi_1(R)$ takes one generator to $t$ and the other to the identity.  The symplectic form can be perturbed slightly so that one or both of the  tori $T_i$ become symplectic. \qed
\end{thm}

\section{Modifying $R$ by Luttinger surgery}\label{sec4}

We can use  Luttinger surgery \cite{Lut} to modify $R$ in a neighborhood of the Lagrangian torus $T_1$ (and $T_2$).

This symplectic construction is carefully explained in \cite{ADK}.  We recall the relevant details for the convenience of the reader.

 An oriented Lagrangian torus $T$ in a symplectic 4-manifold $M$  has a tubular neighborhood symplectomorphic to a neighborhood of  the zero section in its cotangent bundle by the Darboux theorem. Thus if $x,y$ are oriented coordinates on $T$ (i.e. we fix a universal covering $\RR^2\to T$) then $dx,dy$ trivialize the cotangent bundle of $T$, and thereby one obtains a framing $T^2\times D^2\to M$ of a tubular neighborhood of $T$ in $M$ called the {\em Lagrangian framing}.
 As observed in \cite{FS}, this framing can be described by the condition that $T^2\times \{x\}\subset M$ is Lagrangian for
all $x\in D^2$.

Luttinger showed that the manifold obtained by removing the tubular neighborhood $T\times D^2$  from $M$ and regluing using an appropriate  orientation preserving  diffeomorphism $\psi:T\times S^1\to T\times S^1$  yields a new manifold which admits a symplectic structure which agrees with the given symplectic structure on $M-T\times D^2$.  To be precise,  if we denote the generators of $H_1(T)$ (as well as their push offs into $H_1(T\times S^1)$  using the Lagrangian framing by $\alpha, \beta$, and the meridian of $T$ by $\mu$, Luttinger surgery yields a symplectic manifold if
$$\psi_*(\alpha)=\alpha, \ \psi_*(\beta)=\beta, \text{ and } \psi_*(\mu)=a\alpha +b\beta+ \mu.$$  Let  $k=\gcd(a,b)$ and  set $\gamma$ to be the embedded curve $\gamma=\tfrac{1}{k}(a\alpha+b\beta)$, we   denote the resulting symplectic manifold by $M(\gamma,k)$.  Notice that (identifying $\pi_1(T\times S^1)$ with $H_1(T\times S^1)$ and writing multiplicatively)
$$\pi_1(M(\gamma,k))=\pi_1(M-T)/N(\gamma^k\mu).$$
This is because gluing in $T\times D^2$ can be accomplished by gluing in  one 2 handle attached along $\gamma^k\mu$, then gluing in  two 3-handles  and one 4-handle.

 The manifold $M(\gamma,k)$  is  called {\em $1/k$ Luttinger surgery on $T$ along $\gamma$} (the terminology comes from the observation that locally $M$ is obtained by doing a $1/k$ Dehn surgery  and crossing with $S^1$). Notice  that $e(M(\gamma,k))=e(M)$, and Novkov additivity shows that $\sigma(M(\gamma,k))=\sigma(M)$.

\medskip

We apply Luttinger surgery to $T_1\subset R$. Take $\gamma=s$ and $k=1$, and denote the resulting manifold by $ R(s,1)$ by $P$.  Then
$$\pi_1(P)= \pi_1(R-T_1)/N(s\mu_1)= \ZZ t,$$
and $e(P)=8$ and $\sigma(P)=-4$.
Similarly $$\pi_1(P-T_2)=\pi_1(R-(T_1\cup T_2))/ N( s\mu_1 )= \ZZ t$$
and the meridian $\mu_2$ to $T_2$ in $P-T_2$ is trivial, since it is trivial in $\pi_1(R-(T_1\cup T_2))
$.

\bigskip
\noindent{\bf Remark.}  In the preceding paragraph we did not mention the Lagrangian framing, although the reader may check that there is a curve on $T_1$ whose Lagrangian push off is $s\in \pi_1(R-(T_1\cup T_2))$.  An alternative construction of the manifold $R$ is given in Section \ref{mini} in which it is obvious that the Lagrangian framing is obtained by pushing the curve $K$ off itself in the fiber $G$ (and taking the product with $S^1$).

However, note that  if we take an {\em arbitrary}  framing  to push curves off $T_i$, the resulting push offs will differ from the Lagrangian framing push offs by some power of $\mu_1$.  Since $\mu_1=1$ in $ \pi_1(R-(T_1\cup T_2))$, the fundamental group calculation is the same: $s$ is killed, leaving $t$.

\medskip

We will prove that $P$ is minimal  in Theorem \ref{theyaremini} below.

\medskip

We summarize these facts in the following theorem. To simplify the statement we  denote $T_2$ simply by $T$.

\begin{thm}\label{coolthm} There exists a closed minimal symplectic 4-manifold $P$ with fundamental group $\ZZ$, Euler characteristic $e(P)=8$ and signature $\sigma(P)=-4$ which contains a Lagrangian (or symplectic) homologically essential torus $T$ with trivial normal bundle such that \begin{enumerate}
\item  The induced homomorphism $\pi_1(P-T)\to \pi_1(P)$ is an isomorphism,
\item The meridian of $T$ is trivial in $\pi_1(P-T)$, and
\item the induced homomorphism $\pi_1(T)\to \pi_1(P)$ takes one symplectic generator to the generator of $\pi_1(P)$ and maps the other generator to the identity.
 \end{enumerate}\qed
\end{thm}

\bigskip

Our  interest in this manifold is two fold. First, it is the smallest known (to us) symplectic 4-manifold with fundamental group $\ZZ$, where  we measure the size using the Euler characteristic (or equivalently the second Betti number).  (Constructions of symplectic manifolds with fundamental group $\ZZ$ can be found in the literature, e.g. \cite{OS}, \cite{Gompf}, \cite{smith}.) The other reason is that it can be used as a smaller replacement  for the elliptic surface $E(1)$ typically used to control fundamental groups of symplectic 4-manifolds. We will illustrate this in subsequent  sections, using  the following theorem.

\begin{thm} \label{thm4} Let $M$ be a symplectic 4-manifold containing a symplectic torus $T'$ with trivial normal bundle such that $x,y\in \pi_1(M)$ represent the images of the two generators of $\pi_1(T')$. Then the symplectic sum of $P$ and $M$ along $T$ and $T'$, $P\#_T M$, admits a symplectic structure (which agrees with that of $P$ and $M$ away from $T,T'$) symplectic and satisfies
$$e(P\#_T M)=e(M)+8,\ \sigma(P\#_T M)=\sigma(M)-4,\ \text{ and } \ \pi_1(P\#_T M)=\pi_1(M)/N(x)$$
where $N(x)$ denotes the normal subgroup of $\pi_1(M)$ generated by $x$.\qed
\end{thm}

 \bigskip

In a different direction, we can perform Luttinger surgery on both Lagrangian tori $T_1$ and $T_2$.
Suppose we are given   embedded curves $\gamma_1\subset  T_1$ representing
$p_1y_1+q_1s $ and $\gamma_2\subset T_2$ representing $p_2y_2+q_2t $, where $p_1,q_1$ are relatively prime integers and $p_2,q_2$ are relatively prime integers. Suppose further that a pair of integers $k_1,k_2$ are given.

Let $Q=Q(k_1,\gamma_1;k_2,\gamma_2)$ denote the symplectic 4-manifold obtained by performing Luttinger surgery so that $\gamma_1^{k_1}\mu_1 $  and $\gamma_2^{k_2}\mu_2 $ bound  discs after regluing.  We will prove that $Q(k_1,s;k_2,py_2+t)$ is minimal  in Theorem \ref{theyaremini} below. 

  Using Freedman's theorem \cite{Freedman}  and  Taubes's  theorem \cite{taubes}   that minimal symplectic manifolds do not contain $-1$ spheres we immediately conclude the following.

 \begin{thm} \label{cyclic}  The   symplectic manifold $Q(k_1,p_1y_1+q_1s;k_2,p_2y_2 + q_2 t)$
   has Euler characteristic $8$ and signature $-4$, and fundamental group $\ZZ/(k_1q_1)\oplus \ZZ/(k_2q_2)$.   It is minimal when $p_1=0, q_1=1$ and $q_2=1$.

   In particular,  $Q(1,p_1y_1 +s;1, p_2y_2+t)$
 is simply connected
 and hence is homeomorphic to $\CC P^2 \# 5 \overline{\CC P}^2$.  The manifolds $Q(1,s;1,py_2+t)$ (indexed by $p\in \ZZ$) are minimal, hence not diffeomorphic to $\CC P^2 \# 5 \overline{\CC P}^2$.\qed
\end{thm}

Looking from the inside out, we can define $M=Q(1,s;1,t)$. Then $M$ is a minimal symplectic manifold homeomorphic to $\CC P^2 \# 5 \overline{\CC P}^2$ which contains a  (nullhomologous) Lagrangian torus  $T_2$.  Moreover, for each integer $p$, the manifold
$M(p)=Q(1, s; 1, p y_2+t)$ is another minimal symplectic manifold homeomorphic to $\CC P^2 \# 5 \overline{\CC P}^2$, obtained from $M$ by performing Luttinger surgery on $T_2$.
It seems reasonable to conjecture that these are all smoothly distinct.  

\smallskip

The construction of  manifold homeomorphic but not diffeomorphic to $\CP^2\# k\overline{\CP}^2$s for $k\leq 9$ began with Donaldson's seminal example \cite{D} that the Dolgachev surface $E(1)_{2,3}$ is   not diffeomorphic to $ \CP^2\# 9\overline{\CP}^2$. In 1989,  Kotschick \cite{Kot} proved that the Barlow surface is homeomorphic but not diffeomorphic to $\CP^2\# 8\overline{\CP}^2$.  In  2004  J. Park \cite{park1} constructed the first exotic smooth structure on $\CP^2\# 7\overline{\CP}^2$.  Since then Park's results have been expanded upon in \cite{OzS, SS, FS3, PSS}, producing infinite families of smooth 4--manifolds homeomorphic but not diffeomorphic to $\CP^2\# k\overline{\CP}^2$ for $k=5,6,7,8$.  The $n=5$ examples  are not symplectic.

In \cite{A},  Akhmedov  produced the first example of an  symplectic 4-manifold homeomorphic to but not diffeomorphic to  $\CP^2\# 5\overline{\CP}^2$.   He uses a different construction from ours, although the  manifold $W$ is a key ingredient in his construction, as in ours. Very recently Akhmedov, Baldridge,  and Park \cite{AP,bald2} have constructed minimal symplectic manifolds homeomorphic to $\CP^2\# 3\bCP^2$.

There is also an interesting literature on the construction of   manifolds homeomorphic  but not diffeomorphic to  $3\CP^2\#n\bCP^2$ for small $n$.  D. Park \cite{dpark} has constructed many such manifolds and Akhmedov \cite{A} article contains the construction of  a minimal symplectic manifold homeomorphic to $3\CP^2\#7\bCP^2$.  In their recent work  Akhmedov, Baldridge,  and Park \cite{AP,bald2} have constructed minimal symplectic manifolds homeomorphic to $3\CP^2\# 5\bCP^2$.

\bigskip

  Note that the smallest previously known symplectic 4-manifolds with cyclic fundamental group have $e\ge 10$.  There are smooth manifolds with cyclic fundamental group and $e=2$, and any closed 4-manifold with cyclic fundamental group has $e\ge 2$.    In any case, any symplectic 4-manifold with   cyclic fundamental group has $e+\sigma\ge 4$.

\section{minimality}\label{mini}

In this section we prove that the manifolds constructed above are minimal. The key result is the following
theorem.

\begin{thm} \label{min1}  Let $T^4=(S^1\times S^1)\times (S^1\times S^1)$ be given the product symplectic form.  Let $T_1=S^1\times \{a\} \times S^1\times \{a\}$, and $T_2=S^1\times \{b\}\times \{b\}\times S^1$   where $a\ne b$. Denote by $x,y,z,w\in \pi_1(T^4)$ the generators given by the coordinate circles.

Then the symplectic manifold $M$ obtained from  ${1}/{k_1}$ Luttinger surgery on $T_1$ along $x$ and ${1}/{k_2}$ Luttinger surgery on $T_2$ along $w$ is minimal.

\end{thm}
\begin{proof} Notice first that the Lagrangian framings for the tori $T_1$ and $T_2$ are obvious: the parallel Lagrangian tori are obtained by varying the points $a$ and $b$.

Let $M$ denote the symplectic manifold obtained by  ${1}/{k_1}$ Luttinger surgery on $T_1$ along $x$ and ${1}/{k_2}$ Luttinger surgery on $T_2$ along $w$.

It is straightforward to check that $1/k_1$ surgery on $T_1$ along $x$ transforms $T^4=(T^3)\times S^1$ into $N\times S^1$, where $N$ is the 3-manifold that fibers over $S^1$ with monodromy $D_x^{k_1}:T^2\to T^2$, where $D_x$ denotes the Dehn twist on $T^2$ along the first coordinate.  This is explained in \cite[pg. 189]{ADK}.   

View $N\times S^1$ as a trivial circle bundle over $N$.  Then the fiber is represented by the curve $w$ and $T_2$ is the restriction of this circle bundle to $x\times\{(b,b)\}$ in $N$. Luttinger surgery on $T_2$ along $w$ corresponds to twisting the bundle over $x$, i.e. replacing $N\times S^1\to N$ by the $S^1$ bundle $M\to N$ with first Chern class  equal to $k_2$ times the Poincar\'e dual of $x$ in $N$. Details can be found in \cite{B1}. In any case one can check directly from the formula that $M$ has a free circle action which coincides with the action on $N\times S^1$ away from $T_2$.

Thus $M$ is an $S^1$ bundle over a fibered 3-manifold $N$ with fiber a torus. It follows from the long exact sequence of homotopy groups that $\pi_2(M)=0$,  and hence $M$ contains no essential 2-spheres. In particular, $M$ is minimal. \end{proof}

 If  $S:T^4\to T^4$ is a symplectic diffeomorphism then  obviously the manifold obtained by $1/k_1$ surgery on $S(T_1)$ along $S(x)$ and $1/k_2$ surgery on $S(T_2)$ along $S(w)$ is also minimal.    So for example, the linear transformation $\tilde{S}:\RR^4\to \RR^4$ given by $\tilde{S}(x,y,z,w)=(x,y+pz,z, px+w)$  descends to a symplectic diffeomorphism of $T^4$ which leaves $T^1$ and $T^2$ invariant. Thus the manifold obtained by 
 $1/k_1$ surgery on $T_1$ along $x$ and $1/k_2$ surgery on $T_2$ along $x^pw$ is also minimal.

\bigskip

Given an integer $k_1$, let $M_1(k_1)$ denote the symplectic manifold obtained from $T^4$ by two Luttinger surgeries: $1/1$ surgery  on $T_0=S^1\times \{a\} \times \{a\} \times S^1$ along $x$ and $1/k_1$ surgery on $T_1=\{b\}\times S^1\times \{b\}\times S^1$ along $w$.  After suitably permuting coordinates, Theorem \ref{min1} implies that $M_1(k_1)$ is minimal.  Note that $T^4$ contains a pair of  (intersecting) symplectic tori $R_1=\{c\}\times\{c\} \times S^1\times S^1$  and 
$G_1=S^1\times S^1\times \{d\}\times\{d\} $ disjoint from $T_0$ and $T_1$; these tori remain  symplectic in $M_1(k_1)$.

Given a pair of  integers $p$ and $k_2$, let $M_2(p,k_2)$ denote the symplectic manifold obtained from $T^4$ by two Luttinger surgeries: $1/1$ surgery  on $T_3=\{a\}\times S^1\times \{a\}  \times S^1$ along $y$ and $1/k_2$ surgery on $T_2=\{b\}\times S^1\times S^1\times \{b\} $ along $y^pz$.    Theorem \ref{min1} implies that $M_2(p,k_2)$ is minimal.  Again, $T^4$ contains a pair of  symplectic tori $R_2=\{c\}\times\{c\} \times S^1\times S^1$  and 
$G_2=S^1\times S^1\times \{d\}\times\{d\} $ disjoint from $T_3$ and $T_2$; these tori remain  symplectic in $M_2(p, k_2)$.

\bigskip

The fiber sum of $M_1(k_1)$ and $M_2(p,k_2)$ along $R_1$ and $R_2$ is therefore a minimal symplectic manifold by Usher's theorem \cite{usher}.  Denote this manifold by $M(k_1,p,k_2)$. Note that the two symplectic surfaces $G_1$ and $G_2$ can be matched up to produce a symplectic genus 2 surface $G\subset   M(k_1,p,k_2)$ \cite{Gompf}.    Taking the symplectic sum of $M(k_1,p,k_2)$ with the manifold $W$ along the genus two surfaces $G$ and $F$ again produces a minimal symplectic manifold by Usher's theorem since the exceptional spheres in $W$ all intersect the genus two surface $F$. 

A straightforward  examination of the constructions shows that this fiber sum $M(k_1,p,k_2)\#_G W$   
is exactly the symplectic manifold denoted by 
$Q(k_1,\gamma_1;k_2,\gamma_2)$ in Section \ref{sec4}, where $\gamma_1=s$ (so $p_1=0$, $q_1=1$ in the notation of Section \ref{sec4}) and $\gamma_2=py_2 + t$ (so $p_2=p,q_2=1$).   
In fact, note that the fiber sum of  two copies of $T^4$ along $R_1$ and $R_2$ gives $G\times T^2$. The Luttinger surgeries along the tori $T_0$ and $T_3$ correspond exactly to replacing the trivial bundle $G\times S^1\to S^1$ by the fiber bundle with monodromy the Dehn twists $D_1$ and $D_2$. The remaining tori $T_1$ and $T_2$ coincide with the tori labelled $T_1$ and $T_2$ in Section \ref{sec3} and one checks that the surgery coefficients are correct.  We summarize these facts in the following theorem.

\begin{thm}\label{theyaremini} For any triple of integers $k_1,k_2,p$, the symplectic manifold $Q(k_1,s;k_2,py_2+t)$ is  minimal. Its fundamental group is given by 
$$\pi_1(Q(k_1,s;k_2,py_2+t))= \ZZ/k_1 \oplus \ZZ/k_2.$$
It has Euler characteristic equal to $8$ and signature $-4$. 

In particular, the symplectic manifold $R=Q(0,s,0,t)$ with fundamental group $\ZZ^2$ of Theorem \ref{themanifold}  is minimal, the manifold $P=Q(1,s,0,t)$  with infinite cyclic fundamental group  of Theorem \ref{coolthm} is minimal and    the   simply connected  manifolds $M(p)= Q(1,s;1,py_2+t)$ of Theorem \ref{cyclic}  are minimal (and hence exotic) manifolds.  
\qed \end{thm}

 \medskip
 
 \noindent{\bf Remark.} We have not tried to be as general as possible in Theorem \ref{theyaremini}. Presumably a more careful examination will reveal that $Q(k_1,\gamma_1;k_2,\gamma_2)$ is minimal for arbitrary $\gamma_1 $ and $\gamma_2$.    We believe our approach  has much potential: one can consider more complicated monodromies, symplectic sums of more than two copies of $T^4$, and hence surfaces of higher genus and surgeries along more Lagrangian tori, and the use of other manifolds besides $W$.  

\section{More applications.}

\subsection{Free abelian groups}
 We first give an easy application of Theorem \ref{thm4}. Let $M$ denote the 4-torus with its natural symplectic structure and $T'\subset M$ one of the coordinate symplectic tori.

 \begin{cor} $P\#_T M$ is a symplectic 4-manifold with fundamental group $\ZZ^3$, Euler characteristic $8$, and signature $-4$.\qed
 \end{cor}
 The smallest previously known example of a symplectic 4-manifold with fundamental group $\ZZ^3$ has $e=12$ and any such symplectic manifold must have $e\ge 3$ (\cite{BK}).

 \medskip

 More generally, the technique of  \cite[Theorem 20]{BK}  allows us to improve  the construction  of symplectic 4-manifolds with odd rank free abelian groups by taking the fiber sum of $Sym^2(F_n)$ ($F_n$ a surface of genus $n$) with the manifold $P$, rather than the larger manifold $K$ of \cite[Lemma 18]{BK}. We refer the interested reader to loc.cit. for details of the proof of the following corollary.

 \begin{cor} There exists a symplectic 4-manifold $M$ with $\pi_1(M)=\ZZ^{2n-1}$ such that $
 e(M)=11-5n +2n^2$ and $\sigma(M)=-3-n$.\qed
\end{cor}

This is an improvement over the upper bound $\min_{\pi_1(M)=\ZZ^{2n-1}}e(M)\leq 15-5n +2n^2$ of
\cite[Theorem 20]{BK}, but still far from the lower bound $6-7n+2n^2\leq \min_{\pi_1(M)=\ZZ^{2n-1}}e(M)$. Better constructions are needed to decrease  the upper bound (or to increase the lower bound).

 \subsection{Arbitrary fundamental groups.}

As explained in Section 4 of \cite{BK},   the existence of the symplectic manifold $P$ and its symplectic torus $T$ of Theorem \ref{thm4} allows us to   improve  (by 30\%) the main result of loc.cit. to the following theorem.

\begin{thm}\label{30percent} Let  $G$ be a finitely presented group that has a   presentation with $g$
generators and $r$ relations. Then
   there exists a symplectic 4-manifold $M$ with $\pi_1M\cong G$,  Euler
characteristic
$e(M)=12+ 8(g+r)$, and signature $\sigma(M)=-8-4(g+r)$.  Moreover, $M$ contains a symplectic torus  that lies in a cusp neighborhood. \qed
   \end{thm}
 \begin{proof}
 In the proof of \cite[Theorem 6]{BK} a symplectic 4-manifold $N$ is constructed whose fundamental group  contains   classes $s,t, \gamma_1,\cdots,\gamma_{r+g} $  so that $G$ is isomorphic to the quotient of $\pi_1(N)$ by these classes. Moreover, $N$ contains symplectic tori $T_0$, $T_1,\cdots, T_{g+r}$ so that the two generators of $\pi_1(T_0)$ represent $s$ and $t$,  and for $i\ge 1$ the two generators of $\pi_1(T_i)$ represent $s$ and $\gamma_i$.

 Let $E(1)$ denote the elliptic surface with $e(E(1))=12$, the fibration chosen to have (at least) two cusp fibers. let $F_1$, $F_2$ be regular fibers with $F_1$ near one cusp and $F_2$ near another.
Take the fiber sum of $N$ with $E(1)$ along $F_1$, and  $g+r$ copies of the manifold $P$ of Theorem \ref{coolthm} along the tori $T_i$, $i\ge 1$ in such a way that the $\gamma_i$ are killed.  Then a repeated application of Theorem \ref{thm4} computes
$$\pi_1(M)=\pi_1(N)/N(s,t,\gamma_i)=G,\ e(M)=12 +8(g+r), \ \sigma(M)=-8-4(g+r).$$
The torus $F_2$ survives in $M$ as a symplectic torus which lies in a cusp neighborhood.
\end{proof}

 An examination of the proof of Theorem 6 of \cite{BK} shows that Theorem \ref{30percent} can be improved for certain presentations, namely, one can find $M$ so that $e(M)=12 + 8(g'+r)$ and $\sigma(M)=-8-4(g'+r)$, where $g'$ is the number of generators which appear in some relation with negative exponent. Thus if $G$ has a presentation with $r$ relations in which every generator appears only with positive exponent in each relation, then there exists  a symplectic $M$ with $\pi_1(M)=G$, and $e(M)=12+8r,\ \sigma(M)=-8-4r$.  Moreover, using Usher's theorem \cite{usher} one sees that
the manifolds constructed are minimal.

\smallskip

These manifolds can all be assumed to contain a symplectic torus in a cusp neighborhood, since we use one copy of $E(1)$ to get the construction started. Thus the geography results of J. Park \cite{park} can be improved to find a larger region of the $(c_1^2,\chi_h)$ plane for which to each pair of integers in that region one can find   infinitely many non-diffeomorphic, homeomorphic minimal symplectic manifolds with fundamental group $G$.

 \smallskip

 Another (decidedly minor) improvement concerns groups of the form $G\times\ZZ$: for a presentation of $G$ as above there exits a symplectic manifold $M$ with $\pi_1(M)=G\times\ZZ$, $e(M)=  8(g'+r+1)$, and $\sigma(M)=-4(g'+r+1)$. The reason is that one step  in the proof of \cite[Theorem 6]{BK} consists of taking a symplectic sum with $E(1)$ to kill two generators ($t$ and $s$). But to get $G\times \ZZ$ it suffices to kill $t$, for which the manifold $P$ of Theorem \ref{thm4}  can be used instead of $E(1)$. (Of course, if $G$ is the fundamental group of a 3-manifold $Y$ that fibers over $S^1$, then $Y\times S^1$ is a symplectic 4-manifold with fundamental group $G\times \ZZ$ and Euler characteristic zero.)

\smallskip

Suppose $M$ is a symplectic 4-manifold containing  a symplectic torus $T$ with trivial normal bundle so that $\pi_1(M-T)=1$ or $\ZZ$ and so that pushing $T$ into $M-T$ induces a surjection $\pi_1(T)\to \pi_1(M-T)$.  It is easy to see that such an $M$ must have $b^+>1$ and $b^->0$, and hence if $M$ is simply connected $e(M)\ge 6$. If $\pi_1(M)=\ZZ$ such an $M$ must have $e(M)\ge 3$.   The manifold $P$ has $e(P)=8$. Further improvements in the geography problem for symplectic manifolds will be obtained if such an $M$ is found with $e(M)<8$.  The search for such a manifold is a promising direction for future study.

\bigskip
\subsection{Other small simply connected manifolds}

Our manifolds can be used to produce more interesting examples of small simply connected symplectic manifolds.   We consider the case when $b^+=3$. 

  As a warm up, consider the fiber sum of two copies of $P$ along $T$ (see Theorem \ref{coolthm}) If we form the fiber sum along a symplectic diffeomorphism that interchanges (up to sign) the two generators of $\pi_1(T)$,  The resulting manifold is a minimal simply connected symplectic  manifold  homeomorphic to $3\CP^2\# 11\bCP^2$.  One can glue using other maps to obtain manifolds with cyclic fundamental group with $e=16$ and $\sigma=-8$, reminicent of the way  one constructs lens spaces from two solid tori.   In particular there are infinitely many (isotopy classes of) gluings which produce simply connected examples. We do not know if these give distinct diffeomorphism types.

For a smaller example,  we start with a symplectic genus 2 surface $H$ of square zero in 
$T^4\#2\bCP^2$ constructed by symplectically resolving $T^2\times\{(a,a)\}\cup\{(a,a)\}\times T^2$  and blowing up twice.  Thus the homomorphism induced on fundamental groups  by the inclusion $H\subset  T^4\#2\bCP^2$  takes the four standard generators of the genus 2 surface group to the four coordinate generators of $\pi_1(T^4\#2\bCP^2)=\ZZ^4$.  Since $H$ intersects each exceptional sphere  in one point, the meridian of $H$ in $T^4\#2\bCP^2$ is nullhomotopic. Hence an argument exactly like the proof of Lemma \ref{lem1} shows that the symplectic sum of a  manifold $X$ with  $T^4\#2\bCP^2$ along a square zero genus 2 surface has fundamental group obtained from that of $X$ by requiring the four generators of the genus 2 surface in $X$ to commute. 

Consider a parallel copy $F'$ of the symplectic surface $F$ in the manifold $W= (T^2\times S^2)\#4\bCP^2$.    Then the sum $W'$ of $W$ and $T^4\#2\bCP^2$ along $F'$ and $H$ is symplectic and minimal, using Usher's theorem. Moreover, $\pi_1((W-F)\#_{F'}(T^4\#2\bCP^2))=\ZZ^2$, and a quick check shows that  the statement of Lemma \ref{lem1} holds with $W'$ replacing $W$. The difference between the two is that $e(W')=10$ and $\sigma(W')=-6$.  

Taking the fiber sum $R'$ of $W'$ (rather than $W$) with $B$, where $B$ is the manifold constructed in Section \ref{sec3},  and performing Luttinger surgery gives another family $Q'(k_1,\gamma_1;k_2,\gamma_2)$ with $e=14$ and $\sigma=-6$.   As in Theorem 
\ref{cyclic} the manifolds $M'(p)=Q'(1,s,py_2+t)$ are a family of minimal (hence exotic) simply connected symplectic manifolds homeomorphic to $3\CP^2\#9\bCP^2$. We remark that $M'(p)$ is obtained from $M'(0)$ by $1/p$ Luttinger surgery.

\subsection{Some comments on  Seiberg-Witten calculations}

Fix $k$ and consider the  family 
$M(k,p)=Q(k_1,s;1, py_2+t)$ of minimal symplectic manifolds with cyclic fundamental group $\ZZ/k$.   Then $M(k,p)$ can also be described as $1/p$ Luttinger surgery on $C(k,0)$ along $y_2$.  The usual way (\cite{FS2}) to distinguish such families of manifolds is by their Seiberg-Witten invariants, as follows.

 The torus $T_2$ is nullhomologous in $M(k,p)$ for all $p$   and so the Morgan-Mrowka-Szabo formula \cite{MMS} for the Seiberg-Witten invariants of $M(k,p)$  (see also \cite[Theorem 5.3]{FS2})
states that 
$$SW_{M(k,p)}(\kappa_p)= SW_{M(k,0)}(\kappa_0) + p\sum _i SW_{X(k)}(\kappa_X+ i[T])$$
where $X(k)$ is the manifold obtained from  $M(k,0)-nbd(T_2)$ by gluing on $T^2\times D^2$ in such a way that the curve $y_2$ bounds the meridian, and $\kappa_p,\kappa_X$ are  homology classes that correspond (i.e. agree in $H_2(M(k,0),T_2)$). 

However, the manifold $X(k)$ is not obtained by Luttinger surgery, and need not be symplectic. Hence we do not know if $\sum _i SW_{X(k)}(\kappa_X+ i[T])$ is non-zero and it does not immediately follow that the manifolds $M(k,p)$ are distinct.  In fact, when $k=1$ (the case of exotic $\CP^2\#5\bCP^2$) a theorem of Stipcicz and Szabo \cite[Proposition 4.3]{SS} implies  that  $\sum _i SW_{X(1)}(\kappa_X+ i[T])=0$.  In particular, these Seiberg-Witten invariants cannot distinguish the manifolds $M(1,p)$.

 \vfill\eject




\end{document}